\documentclass[12pt]{amsart}
\usepackage{amsmath,amsthm}
\usepackage{mathrsfs}
\input amssym.def
\input amssym
\usepackage{enumerate}
\usepackage{hyperref}
\hypersetup{
    pdftitle={Super},    
    pdfauthor={G\"otz Pfander}{Ilya Krishtal},     
    pdfsubject={Superresolution},   
    pdfcreator={Ilya Krishtal},   
    pdfproducer={Ilya Krishtal}, 
    pdfkeywords={Banach modules}{Beurling spectrum}{Wiener's lemma}, 
     colorlinks,%
    citecolor=blue,%
    filecolor=black,%
    linkcolor=magenta,%
    urlcolor=black
}

\usepackage{algorithm}
\usepackage{algorithmicx}
\usepackage[noend]{algpseudocode}
\usepackage{tikz}
\allowdisplaybreaks

\chardef\bslash=`\\ 

\makeatletter
\def\verbatim{\interlinepenalty\@M \@verbatim
   \leftskip\@totalleftmargin\advance\leftskip2pc
   \frenchspacing\@vobeyspaces \@xverbatim}
\makeatother
\newtheorem{thm}{Theorem}[section]
\newtheorem{cor}[thm]{Corollary}
\newtheorem{lem}[thm]{Lemma}
\newtheorem{prop}[thm]{Proposition}

\theoremstyle{definition}
\newtheorem{defn}{Definition}[section]

\theoremstyle{remark}
\newtheorem{rem}{Remark}[section]
\newtheorem{exmp}{Example}[section]

\numberwithin{equation}{section}

\newcommand{\begeq}{\begin {equation}}
\newcommand{\eq}{\end{equation}}
\newcommand{\bs}{\begin {split}}
\newcommand{\es}{\end{split}}
\newcommand{\bp}{\begin {prop}}
\newcommand{\ep}{\end {prop}}
\newcommand{\bt}{\begin {thm}}
\newcommand{\et}{\end {thm}}
\newcommand{\bc}{\begin {cor}}
\newcommand{\ec}{\end {cor}}
\newcommand{\bl}{\begin {lem}}
\newcommand{\el}{\end {lem}}
\newcommand{\bpf}{\begin {proof}}
\newcommand{\epf}{\end {proof}}
\newcommand{\bi}{\begin {itemize}}
\newcommand{\ei}{\end {itemize}}
\newcommand{\ben}{\begin {enumerate}}
\newcommand{\een}{\end {enumerate}}
\newcommand{\brem}{\begin {rem}}
\newcommand{\erem}{\end {rem}}

\newcommand{\bd}{\begin {defn}}
\newcommand{\ed}{\end {defn}}

\newcommand{\bex}{\begin {exmp}}
\newcommand{\eex}{\end {exmp}}

\newcommand{\norm}[1]{\left\|{#1} \right\|}

\newcommand{\la}{\langle}
\newcommand{\ra}{\rangle}

\newcommand{\fhat}{\hat{f}}

\newcommand{\A}{\mathcal A}

\newcommand{\F}{\mathcal{F}}
\newcommand{\I}{{\mathcal I}}

\newcommand{\Y}{\mathcal{Y}}
\newcommand{\TTT}{{\T\kern-.44em \T}}
\newcommand{\tTTT}{\widetilde{\T\kern-.44em \T}}

\newcommand{\KK}{\mathcal K}

\newcommand{\ZZ}{{\mathbb Z}}
\newcommand{\TT}{{\mathbb T}}
\newcommand{\RR}{{\mathbb R}}
\newcommand{\CC}{{\mathbb C}}
\newcommand{\NN}{{\mathbb N}}

\newcommand{\X}{\mathcal{X}}

\newcommand{\GG}{\mathbb{G}}

\newcommand{\G}{\mathbb G }

\newcommand{\HH}{{\mathcal H}}

\newcommand{\s}{\sigma}
\renewcommand{\l}{\lambda}
\renewcommand{\L}{\Lambda}
\renewcommand{\a}{\alpha}

\newcommand{\g}{\gamma}
\newcommand{\Gg}{\widehat{\G}}
\newcommand{\T}{\mathcal{T}}

\DeclareMathOperator{\supp}{supp}

\newcommand{\aee}{\mbox{\rm \ae}}

\begin{document}

\title {Prony's Method in Banach Modules}

\author{ Ilya A. Krishtal and G\"otz E. Pfander }
\address{Department of Mathematical Sciences, Northern Illinois University, DeKalb, IL 60115 \\
email: ikrishtal@niu.edu}
\address{Lehrstuhl Wissenschaftliches Rechnen, Mathematisch-Geographische Fakult\"at, Karholische Universit\"at Eichst\"att Ingolstadt, Germany\\
email: goetz@pfander.org}


\date{\today }




\begin{abstract}
We show that the classical Prony's method for recovery of a sparse signal from its consecutive Fourier coefficients  can be viewed as a spectral identification problem for an unknown restriction of a known linear operator. This presents a unified point of view on various existing and novel generalizations and applications of the method, some of which are discussed in this paper.  
\end{abstract}
\maketitle

\section{Introduction}

In 1795, Prony demonstrated in his work \cite{P795} that it is possible to reconstruct any $s$-sparse $d$-dimensional vector from any $2s$ of its consecutive Fourier coefficients. Since then, Prony's ideas have been expanded in various directions, leading to advancements in several sub-fields 
of applied and computational mathematics (see \cite{BDVMC08, CL18, JLM19, PP13, S17} and references therein for a small sample of these developments).

In the paper \cite{AK16}, the authors take Prony's method a step further and present a novel perspective on it. They show that the core of Prony's method can be seen as a spectral identification problem, where one aims to identify an unknown restriction of a known linear operator in a Euclidean space. In this paper, we extend this viewpoint to a much more abstract level, considering Banach modules over group algebras instead of just Euclidean spaces.

By generalizing Prony's method in this way, we achieve a unification of various existing extensions of Prony's method \cite{BDVMC08, PP13} and address other related problems \cite{AK16, JLM19}. One of the problems discussed in the paper is the identification of operators represented as linear combinations of time-frequency shifts \cite{gro01,GP14,Heil2003, Kozek1998}. These operators are commonly used to model time-varying communication channels \cite{Vtr71}. Accurately identifying the channel is crucial for reliable communication, and it typically involves observing a limited number of input and output pairs of the operator \cite{GP14,KP06,PW06,HeckelBoel13,PZ14b}. See \cite{WPK15} for an overview and an history of the subjects, and additional references. In Section \ref{TFS} of the paper, we present a novel approach to tackle this problem that is based on a special case of our abstract method.

To provide the necessary background, we introduce the relevant concepts in the beginning  of Section \ref{set}. The rest of the section is the centerpiece of the paper where our abstract approach is presented. We proceed to demonstrate how Prony's method and some of its generalizations can be derived as special cases in Section \ref{Prony}. Finally, in the concluding remarks, we discuss further potential applications of the method.

\section{Banach modules and spectral identification}\label{set}

We begin by introducing the notation and relevant notions of the spectral theory of Banach modules. We refer to \cite{B04, BK05, BK14} and references therein for more information on the subject.

The symbol $\X$ will denote a complex Banach space and $B(\X)$ will be the Banach algebra of all
bounded linear operators on $\X$. By $\G$ we will denote a locally compact Abelian group (LCA-group)
and by $\Gg$ -- its \emph{Pontryagin dual}, i.e.~the group of all continuous unitary characters of $\G$. 
We will write the operation additively on all LCA-groups, except when the group $\TT = \{\theta\in\CC:\ |\theta| = 1\}$ is used.

\begin{defn} \label{wegt}
A  \emph{(Beurling) weight} is a (Haar) measurable locally bounded even function $w: \G\to [1,\infty)$ such that
\[ w(g_1+g_2)\leq w(g_1)w(g_2),\ \mbox{ for all }\ g_1,g_2\in\G. \]
\end{defn}
 
We denote by $L_w(\G)$ the Banach algebra of (equivalence classes of) complex-valued functions $w$-integrable
on $\G$  with respect to the normalized Haar measure. The norm is given by
\[\norm{f}_w = \int_\G |f(g)|w(g)dg <\infty,\] and the role of
multiplication in $L_w(\G)$ is played by the convolution 
\[(f\ast h) (g)= \int_\G f(t)h(g-t)dt.\]
These algebras are typically called \emph{Beurling algebras}  \cite[and references therein]{F79g,F84,HNR98,K09,L85,O75, R68,RS00,St76}. If $w \equiv 1$, we obtain the standard group algebra, that is, $L_1(\G) \equiv L^1(\G)$.

In this paper, we typically assume that the weight $w$ is \emph{non-quazi-analytic} \cite{D56, L85, LMF73},  i.e.~ 
it satisfies the \emph{Beurling-Domar condition}
\begeq\label{bdc}
\sum_{n=0}^\infty \frac{\ln w(ng)}{1+n^2} < \infty \quad \mbox{for all } g\in\G,\quad ng = \underbrace{g+g+\ldots+ g}_{n\, {\rm times}}.
\eq
This condition ensures  that the algebra $L_w(\G)$ is regular and its spectrum 
 is isomorphic to the dual group $\Gg$ \cite{D56, K09}. 
 
We denote by $\fhat : \Gg \to \CC$ the \emph{Fourier transform} of $f \in L_w(\G)$:
\[\hat f(\g) = \int_{\G} f(g)\g(-g) dg,\ \g\in\Gg.\]
 The \emph{inverse Fourier
transform} of a function  $h: \Gg \to \CC$ is denoted by $\check h$ or $h^\vee$.  

For some results, we may also require that  the weight $w$ satisfies \emph{Ditkin's condition}  \cite{GRS64}, which  states that
given $\g\in\Gg$ and a function $f\in L_w(\G)$ with $\hat f(\gamma) = 0$ there exists a sequence $(f_n)_{n\in\NN}$ of functions in $L_w(\G)$ such that $\hat f_n \equiv 0$ in a neighborhood of $\g$ and
$\lim f_n*f = f$. Ditkin's condition is satisfied \cite{B82} if, for example,  
 \begeq\label{dc}
 \lim_{n\to \infty}\frac{w(ng)}{n} =0\ \mbox{ for all } g\in\G.
 \eq

In this paper, the space $\X$ is endowed with a Banach $L_w(\G)$-module structure associated with some representation
$\T: \G \to B(\X)$. We shall occasionally use the notation $ (\X,\T)$ for such modules. If $\X$ is an appropriate space of functions (or distributions) on $\G$, the following representations are used most often.

The \emph{translation} representation $\T = T: \G\to B(\X)$ is given by
\begeq\label{trans1}
T(t)x(s) = x(s+t),\ x\in\X, \ s,t\in\G,
\eq 
and the \emph{modulation} representation $\T = M: \Gg \to B(\X)$ is defined by
\begeq\label{mod1}
M(\xi)x(s) = \xi(s)x(s),\ x\in\X, \ \xi\in\Gg, \ s\in \G.
\eq 

We will only consider representations $\T: \G \to B(\X)$ that are \emph{strongly continuous}; in other words, we require that the function $x_\T: \G\to\X$ given by $x_\T(g) =\T(g)x$ is continuous for every $x\in\X$. We also assume $ \norm{\T(g)}\le w(g)$ for all $g\in\G$.
The  Banach $L_w(\G)$-module structure \cite{BR75, F84, R67}  on $\X$ is then defined via
\begeq\label{scdef}
fx =  \int_\G  f(g)\T(-g)x dg,\quad f\in L_w(\G), \ x\in\X.
\eq
We will always assume that this structure is \emph{non-degenerate}, that is, $fx = 0$ for all $f\in L_w(\G)$ implies $x=0$. We also note the following important property of the defined module structure:
\begeq\label{redef}
\T(g)(fx) = \left(T(g)f\right)x = f\left(\T(g)x\right),
\eq  
where  $T$ is the translation representation (\ref{trans1}).

\bd
The \emph{Beurling spectrum} of a set $N \subseteq (\X,\T)$ is the subset $\L(N)= \L(N,\T)$ of the dual group $\Gg$ 
given by
\[\{\g\in\Gg: \mbox{ given any } f\in L_w(\GG) \mbox{ with } \hat f(\g) \ne 0 \] \[\mbox{ there is } x\in N \mbox{ such that } fx \neq  0\}.
\]
When $N = \{x\}$ is a singleton, we shall write $\Lambda(x)$ instead of $\L(\{x\})$.
\ed

The notion of the Beurling spectrum serves as a natural generalization of the notion of a support.
Among many things the Beurling spectrum may coincide with, specific choices of Banach modules  yield the support of a function, the support of the (distributional) Fourier transform, the set of indices of non-zero matrix diagonals, etc \cite{BK05}. The following example serves as an illustration.

\begin{exmp}
    Let $\X = L^2(\G)$. If $\T = M$ is the modulation representation \eqref{mod1}, then the $L^1(\Gg)$-module structure of $(\X, M)$ is given by $(fx)(s) = \widehat f(s) x(s)$ so that $\Lambda(x) = \supp x$. If $\T = T$ is the translation representation \eqref{trans1}, then the $L^1(\G)$-module structure of $(\X, T)$ is given by the convolution $(fx)(s) =  (f*x)(s)$. As $f*x=0$ if and only if $\widehat f \,\widehat x=0$ we conclude that  $\Lambda(x) = \supp \widehat x$.  
\end{exmp}

Given a closed set $F\subset \Gg$ we shall denote by $\X(F)$ the (closed) \emph{spectral submodule} of
all vectors $x\in\X$ such that $\L(x)\subseteq F$. 
 If $F = \{\g\}$ is a singleton, we shall write $\X_\g$ instead of $\X(\{\g\})$. 
 
 \brem\label{dr} A submodule $\X_\g$, $\g\in\Gg$,  contains all $x\in\X$ such that $\T(g)x = \g(g)x$ for all $g\in\G$. Moreover, if Ditkin's condition is satisfied,  $\X_\g$ does not contain any other vectors (see, e.g., \cite{B78}).
 \erem
 
The key objects of this paper are the spectral submodules $\X(F)$ generated by a \emph{finite} set $F$ of a fixed cardinality $\ae$.

Given some $x\in \X(F)$, the goal is to use certain ``measurements'' of $x$ to identify initially the set $F$ and ultimately the vector $x$ itself. Only a finite number of measurements shall be used. Therefore, it is natural to assume that $F\subseteq \Omega$ for some known subset $\Omega\subseteq \Gg$, and each $\X_\g$, $\g\in\Omega$, has a finite dimension $m_\g \le M\in\NN$. The set $\{x_\g^m:\ m =1,\ldots, m_\g\}$ will be a (known) basis in $\X_\g$, $\g\in\Gg$. Thus, our signal model is
\begeq\label{signal}
x = \sum_{\g\in F}\sum_{m=1}^{m_\g} c_{\g m} x_\g^m,\ c_{\g m}\in \CC.
\eq
The set $F$ and the coefficients $c_{\gamma m}$ are unknown and to be recovered from certain ``measurements'' of $x$.
It will  be assumed that for each $\g\in F$ there exists $m\in\{1,\ldots, m_\g\}$ such that $c_{\g m}\neq 0$.

To compute the measurements, we shall choose two linear operators $A:\X\to \CC^S$, $S\in\NN$, and $B\in B(\X)$. The measurements $y_\ell$ shall be of the form
\begeq\label{meas}
y_\ell = AB^\ell x, \ \ell =0,1,\ldots, L.
\eq
The choice of the operator $A$ is predicated by two goals it needs to achieve. 
One of them is to allow for unique
 identification of the coefficients $c_{\gamma m}$ once the set $F$ is recovered. In general, one can only ensure this by requiring that all matrices in a certain class have a left inverse. In most of the specific examples we present in this paper, however, we check that this condition is satisfied. The other goal in the choice of the operator $A$ is to make sure that it does not obscure the recovery of $F$ made possible
 via the design of the operator $B$. Sufficient conditions for achieving this goal are included in the statements of the general results.

The operator $B$ shall have the form
\begeq\label{opB}
B = \sum_{n=1}^N b_n \T(g_n), \ g_n\in\G, b_n\in\CC\setminus\{0\}. 
\eq
By $B_F$ we shall denote the restriction of the operator $B$ to $\X(F)$. By definition,  all spectral submodules are invariant for $B$ of the form \eqref{opB}, so that the restrictions are well defined. The following important result, which is a special case of the spectral mapping theorem, indicates why using $B$ of the form \eqref{opB} may lead to the recovery of $F$.

\bt[{\cite[Theorem 3.3.14]{B04}}]\label{spec1}
The spectrum $\s(B_F)$ of the operator $B_F$ satisfies
\begeq\label{spec2}
\s(B_F) =\overline{\left\{\sum_{n=1}^N b_n \g(g_n): \g\in F\right\}}.
\eq
\et
The right-hand-side of \eqref{spec2} motivates us to introduce the function $h =h_B: \Gg\to\CC$ given by
\begeq\label{funh}
h(\g) = \sum_{n=1}^N b_n \g(g_n), \ \g\in\Gg.
\eq
With this notation, using the fact that the image of a compact set $F$ under a continuous function $h$ is compact, we can write \eqref{spec2} more succinctly:
\begeq\label{spec3}
\s(B_F) 
= h_B(F).
\eq

In view of Theorem \ref{spec1}, we would like to be able to choose the coefficients $b_n$ in \eqref{opB} so that $\s(B_F) = {h(F)}$ completely determines $F$. 
If this is, indeed, possible, the method of \cite{AK16} outlined below can be used to determine the ``observable'' part of $F$. For convenience, we always choose the coefficients $b_n$ in such a way that $0\notin \s(B_F)$ for any nonempty $F\subseteq \Omega$.

Observe that we  have $\dim \X(F) \le\ae M$. Thus, for $x\in \X(F)$, the system of vectors 
$\{x, Bx, \ldots, B^{\ae M}x\}$ is always linearly dependent and we can find the coefficients
$\a_0$, $\a_1$, \ldots, $\a_{\ae M}$ such that $\a_{\ae M} = 1$ and
\begeq\label{msym}
 \sum_{\ell=0}^{\ae M} \a_\ell y_{\ell+k} =\sum_{\ell =0}^{\ae M} \a_\ell AB^{\ell+k} x = 0,\ k = 0,1, \ldots
\eq
We let 
\begeq\label{pmineq}
p_{\min}(z) = \sum_{\ell =0}^{\ae M} \a_\ell z^\ell, \quad z\in\CC,
\end{equation}
where $\a_\ell$ are such that \eqref{msym} holds; if the choice is not unique, we pick $0 = \a_{0} = \a_{1} = \ldots$ for as many coefficients as possible. Observe that $p_{\min}$ can be found from the measurements \eqref{meas} with $L = 2\ae M-1$.

The following result is now an immediate consequence of \cite[Proposition 2.4]{AK16}.

\bp\label{mainprop}
The set $R_{\min}$ of all non-zero roots of the polynomial $p_{\min}$ is a subset of $\s(B_F)$.
\ep

The above proposition can be considerably strengthened when all singleton spectral submodules are one-dimensional, i.e.~$M = 1$. In that case, we have $\X_\g = \{\a x_\g: \a\in\CC\}$ and $\T(g)x_\g = \g(g)x_\g$, $g\in\G$.

\bt\label{easythm}
Assume that the restriction of $h$ to $\Omega$ is one-to-one and $0\notin {h(\Omega)}$. Assume also that $M = 1$ and for each $\g\in\Omega$ we have $Ax_\g \neq 0$. Then $R_{\min} = h(F)$ and $F$ can be recovered from the measurements \eqref{meas} with $L = 2\ae-1$.
\et

\bpf
From \eqref{msym} and  \cite[Proposition 2.4]{AK16}, we have
\begeq\label{defeq0}
 \sum_{\ell=0}^{\ae} \a_\ell AB^{\ell+k}x = \sum_{\g \in F}  c_\g h^k(\g)p_{\min}(h(\g))
 Ax_{\g}= 0
\eq
 for all $k=0, 1,\ldots$ Our assumptions on $h$ ensure that the Vandermonde matrix generated by $\{h(\g): \g\in F\}$ is invertible. Since we have also assumed that $c_\g Ax_\g\neq0$, \eqref{defeq0} may only hold if $p_{\min}(h(\g)) =0$ for all $\g\in F$. Hence, our choice of $p_{\min}$ ensures that $R_{\min} = h(F)$. The set $F$ can now  be recovered since the restriction of $h$ to $\Omega$ is one-to-one.
\epf

To obtain an analogous result for the case  $M > 1$ we will need the following definitions.

\bd
A spectral value $\g\in F$ is called \emph{observable} if the restriction of $h$ to $\Omega$ is one-to-one,  $0\notin {h(\Omega)}$, and $h(\g)\in R_{\min}$. The set $F$ is called \emph{
observable} if all $\g\in F$ are observable.
\ed

\bd
A \emph{Krylov subspace} of order $r$ generated by an operator $B$ and a vector $x\in\X$ is
\[
\mathcal K_r(B, x) = {\rm  span}\{x, Bx, . . . , B_{r-1}x\}. 
\]
The \emph{maximal Krylov subspace} will be denoted by $\mathcal K_{\infty}(B,x)$.
\ed

\bd\label{defminpoly} 
For $x\in\X$, the {\em $(A,B,x)$-annihilator}, denoted by  $p^B_{A,x}$, is the monic polynomial of the smallest degree among all the polynomials $p$ such that $Ap(B)\KK_{\infty} (B,x) =\{0\}.$
\ed

We have the following sufficient criterion for observability.

\bp\label{propobs}
Assume that the restriction of $h$ to $\Omega$ is one-to-one and $0\notin {h(\Omega)}$. Assume also that $\g\in F$ is such that there is a (non-zero) eigenvector $x_\g^1\in \mathcal K_{\infty}(B,x)$ satisfying $B x_\g^1 = h(\g)x_\g^1$ and $Ax_\g^1 \neq 0$. Then $\g$ is observable.
\ep

\bpf
By construction, we have $p_{\min} (z)= z^kp^B_{A,x}(z)$ for some $k \ge 0$. Since $x_\g^1\in \mathcal K_{\infty}(B,x)$, it follows that
\[
0 = Ap_{\min}(B) x_\g^1 = p_{\min}(h(\g))Ax_\g^1 = (h(\g))^k p^B_{A,x}(h(\g))Ax_\g^1.
\]
Since $Ax_\g^1 \neq 0$ and $0\notin {h(\Omega)}$, we must have $p_{\min}(h(\g))=0$. It follows that $h(\g)\in R_{\min}$.
\epf

The following lemma is implied by standard linear algebraic facts. We provide a proof for completeness. 

\bl\label{l4}
For any $\g \in F$ there is $x_\g^1 \in \mathcal K_{\infty}(B,x)\setminus\{0\}$ such that $B x_\g^1 = h(\g)x_\g^1$. 
\el

\bpf
Pick $\g\in F$ and let $E = h(F)\setminus\{h(\g)\}$. We will use a polynomial $q$ defined in the following way. If $E =\emptyset$, we let $q \equiv 1$. Otherwise, we let
\[
q(z) = \prod_{\tau\in E}(z-\tau).
\]
By construction, we have $q(B)x \in \mathcal K_{\infty}(B,x)\setminus\{0\}$.
Let $\Y$ be the submodule generated by $q(B)x$. From Theorem \ref{spec1}, we deduce that the restriction $B_{\Y}$ of $B$ to $\Y$ has the spectrum $\s(B_{\Y}) = \{h(\g)\}$. The conclusion of the lemma now follows from the fundamental theorem of algebra.
\epf

Combining Proposition \ref{propobs} with Lemma \ref{l4} and Remark \ref{dr}, we get

\bt\label{hardthm}
Assume that the restriction of $h$ to $\Omega$ is one-to-one and $0\notin {h(\Omega)}$. 
Then  $\g \in F$ is observable if and only if $Ax\neq 0$ for each $x\in\X$ such that $\T(g)x = \g(g)x$ for all $g\in\G$.
\et

The results of this section yield the following Algorithm \ref{alg:rec}, which is essentially the same as the classical Prony algorithm.

 	\begin{algorithm}[h]
\caption{Spectral recovery.}
	\label{alg:rec}
	\begin{algorithmic}[1]
		\State \textbf{Goal: }Find 
  all elements $x\in\X$ of the form \eqref{signal} that satisfy the given measurements \eqref{meas}.
		\State \textbf{Input: } The measurements \eqref{meas}, function $h$ satisfying conditions of Theorems \ref{easythm} and \ref{hardthm}, maximal cardinality $\aee$ of the spectrum $F$, and maximal dimension $M$ of the spectral submodules $\X_\g$.
    \State {Set up the linear system \eqref{msym} and solve it with $0 = \a_{0} = \a_{1} = \ldots$ for as many coefficients as possible.}
    \State {Form the polynomial $p_{\min}$ given by \eqref{pmineq} and find the set $R_{\min}$  of its non-zero roots.}
    \State {Find $F = h^{-1}(R_{\min})$. }
	\State {Find all possible solutions for $c_{\g m}$ by solving the linear system yielded by \eqref{meas} when $x$ with the recovered $F$ is plugged in.}
	\State	
 \textbf{Output: } All $x\in\X$ of the form \eqref{signal} that satisfy  \eqref{meas}.   
	\end{algorithmic}
\end{algorithm}

\begin{rem}
Algorithm \ref{alg:rec} fully determines the spectrum $F$ of any $x\in\X$ of the form \eqref{signal} satisfying \eqref{meas}. It is also clear, that
 the final linear system in the algorithm is guaranteed to have a solution. In general, however, there is no guarantee that the solution is unique, i.e.~there may be infinitely many $x\in\X$ of the form \eqref{signal} that satisfy the given measurements \eqref{meas} (they will all have the same spectrum $F$). We show in the following section that sufficient conditions for the uniqueness of the solution may often be formulated in terms of an appropriate choice of the operator $A$.
\end{rem}

\section{Reduction to Prony's method.}\label{Prony}

In this section, we recover  
 Prony's method \cite{BDVMC08} and develop its more general versions. 

\subsection{Standard Prony's method.}

(Maybe we should also address the on-grid case, for flow of reading.  We could even make a short section, but I think that would not work.)

Here, we let $\G = \RR$, $\Gg \simeq \RR$,
 $\X = C_{ub}(\RR)$ -- the space of all bounded uniformly continuous functions on $\RR$, and $F \subset [0,1)$. The role of $\g(g)$, $g\in\G$, $\g\in \Gg$, in this case, is played by   $e^{2\pi i g\g}$, $g, \g\in\RR$.
 The signal model is then
\begeq\label{pronysig}
x(t) = \sum_{\g\in F} c_\g e^{2\pi i t\g},\ t\in\RR,
\eq
i.e.~$x$ is an almost periodic function with a finite spectrum.
The representation $\T$ is the translation $\T(t)y(s) = y(s+t)$, $y\in\X$, so that the module structure is the convolution:  $fy = f*y$,  $f\in L^1(\RR)$. It is then easy too see from Remark \ref{dr}  that letting $x_\gamma(t) = e^{2\pi i t\g}$
yields $\X_\gamma = \{a x_\gamma: a\in\CC\}$. Thus, \eqref{pronysig} is indeed a special case of \eqref{signal} and all spectral submodules $\X_\g$ are one-dimensional.

We now let $B = \T(1)$ and $A: \X\to \CC$ be given by $Ay = y(0)$.
This yields  standard Prony's measurements
\begeq\label{pronymeas}
y_\ell   = AB^\ell x= x(\ell), \ \ell = 0, 1, \ldots.
\eq

The function $h: \RR\to \CC$ in \eqref{funh} becomes $h(\g) = e^{2\pi i \g}$. Clearly, $h$ is one-to-one on $\Omega = [0,1)$, and $0\notin h(\Omega) = \TT$. Thus, Theorem \ref{easythm} applies ascertaining the validity of the standard off-grid Prony's method.

To summarize, the system \eqref{msym} in this case becomes
\[
 \sum_{\ell=0}^{\ae} \a_\ell y(\ell+k) = 0, \ k = 0, 1, \ldots, \ae-1.
\]

 The solution with as many $0 = \a_0 = \a_{1} = \ldots$ as possible provides the coefficients of $p_{\min}$. Finding the non-zero roots   of $p_{\min}$, one   obtains the spectrum of $x$ from $R_{\min} = h(F)$. To find $x$ itself it remains to solve another linear system, which is guaranteed to have a unique solution since its matrix is Vandermonde.
 
 \brem
 One could make another choice of $B$ above. Potentially this allows one to better distinguish the frequencies that are close by.   
 \erem

 \subsection{Exponential sums with polynomial coefficients.}
A more general form of the off-grid Prony's method arises if we let
the signal $x$    be a continuous function of the form
\begeq\label{genpronysig}
 x(t) = \sum_{\g\in F} q_\g(t)e^{2\pi i\g t}, 
\eq
 where each $q_\g$ is a (non-zero) polynomial of degree $m_\g\le M$. 
The space $\X$ is then the closed linear span of such signals in
 \[
\Y = \left\{y: \frac y{(1+|\cdot|)^M}\in C_{ub}(\RR)\right\},
 \]
 with  $ \|y\|_{\Y} = \sup_{t\in\RR} |y(t)|{(1+|t|)^{-M}}$. We again have that the representation $\T$ is the translation $\T(t)y(s) = y(s+t)$, $y\in\X$, and the module structure is the convolution:  $fy = f*y$. It is crucial, however, that this time the module is not over $L^1(\RR)$ but only over $L_w(\RR)$, where $w(t) = (1+|t|)^M$. Thus, even though the functions $x_\g = x_\g^1$ given by $x_\g(t) = e^{2\pi i \g t}$ are still the eigenvectors of the representation $\T$  and, therefore, belong to  the spectral submodules $\X_\gamma$, they no longer span them. In this case, the submodules $\X_\gamma$ are not one-dimensional. In fact, one can check that the functions $x_\gamma^m$, $m=1,\ldots M+1$, given by $x_\g^m(t) = t^{m-1}e^{2\pi i \g t}$ form a basis of $\X_\g$. In particular, to see that $x_\g^m \in \X_\g$, notice that integration by parts yields $(f*x_\g^m)(t) = e^{2\pi i \g t} \sum_{k=0}^m t^k\widehat f_k(\g)$, where $\supp \widehat f_k = \supp \widehat f$ for each $k = 0, \ldots, m$. Choosing $f$ with $\g\notin \supp \widehat f$, therefore, yields $f*x_\g^m = 0$ resulting in $\Lambda(x_\g^m) = \{\g\}$.
It follows that \eqref{genpronysig} is also a special case of \eqref{signal}.
 
 We now choose $A$ and $B$ the same way as in the previous subsection. We have that all eigenspaces of $B$ are of the form $\{\a x_\g^1: \a\in\CC\}$. 
 Moreover, $Ax_\g^1 = 1\neq 0$, and Theorem \ref{hardthm} applies yielding the following corollary.
 
 \bc
Let $x(t) = \sum\limits_{\g \in F} q_\g(t)e^{2\pi i \g  t}$, where $t, \g\in \RR$,  $q_\g$ is a non-zero polynomial of degree at most $M$, and $F$ is a finite subset of $[0,1)$ of cardinality $\ae$. 
Then $x$ can be completely recovered from its values $x(\ell)$, $\ell =0,1,\ldots 2\ae (M+1)-1$.
\ec

\bpf
From Theorem \ref{hardthm} we deduce that $F$ is 
observable. It remains to prove that the polynomials $q_\g$ can be recovered from the values of $x$ once $F$ is known. This is shown in the following Lemma \ref{aplem}. 
\epf

\bl\label{aplem}
Let $x(t) = \sum\limits_{n =1}^N q_n(t)(\theta_n)^t$, where $t\in \RR$, $\theta_n\in\TT$ are distinct, and $q_n$ are unknown non-zero polynomials of degree at most $M$,  $n =1, 2, \dots, N$. 
Then $x$ can be completely recovered from its values $x(\ell)$, $\ell =0,1,\ldots N(M+1)-1$.
\el

\bpf
The values $x(\ell)$ allow one to write a linear system of equations for the unknown coefficients of the polynomials $q_n$.
The lemma follows from the invertibility of the matrix of that system, which can be written as 
$
(V_1 \ V_2 \ldots V_N),
$
where $K = N(M+1)-1$ and 
\[
V_n = 
\left(
\begin{array}{ccccc}
1 & 0 & 0 & \dots & 0  \\
\theta_n & \theta_n & \theta_n & \dots &\theta_n  \\
\theta_n^2 & 2\theta_n^2 & 2^2\theta_n^2 & \dots &2^{M}\theta_n^2  \\ 
\vdots& \vdots&\vdots&\ddots&\vdots  \\
\theta_n^K & K\theta_n^K & K^2\theta_n^K & \dots &K^{M}\theta_n^K  
\end{array}
\right), \quad
n= 1, \ldots, N.
\]

By Gaussian elimination in each $V_n$, 
$n=1, \dots, N$, the above matrix reduces to  $(\widetilde V_1 \ \widetilde V_2 \ldots \widetilde V_N),$ where
\[
\widetilde V_n =
\left(
\begin{array}{ccccc}
1 & 0 & 0 & \dots & 0  \\
\theta_n & 1 & 0 & \dots &0  \\
\theta_n^2 & 2\theta_n & 1 & \dots &0  \\ 
\vdots& \vdots&\vdots&\ddots&\vdots \\
\theta_n^M & M\theta_n^{M-1} & {M\choose2}\theta_n^{M-2} & \dots &1  \\ 
\vdots& \vdots&\vdots&\ddots&\vdots \\
\theta_n^K & K\theta_n^{K-1} & {K\choose2}\theta_n^{K-2} & \dots &{K\choose  M}\theta_n^{K-M} 
\end{array}
\right), \quad
n= 1, \ldots, N.
\]
We show that this matrix is invertible by augmenting the standard elementary argument for computing the determinant of a Vandermonde matrix with Gaussian elimination in each $\widetilde V_n$, $n= 1, \ldots, N$. 

The proof is by induction in $M$. In the base case $M=0$ the matrix is Vandermonde and is, therefore, invertible. The beginning of the inductive step is precisely the same as in the Vandermonde induction: the first column of $\widetilde V_1$, i.e.~the block corresponding to $\theta_1$,  is subtracted from the first columns of all other blocks and common factors are eliminated in each modified column; then from each row except the first one we subtract the previous row multiplied by $\theta_1$. After this step, the  $M+1$ columns corresponding to $\theta_1$
become
\[
\left(\!\!
\begin{array}{cccccccccccccccc}
1 & 0 & 0 & \dots & 0 \\
0 & 1 & 0 & \dots &0  \\
0 & \theta_1 & 1 & \dots &0  \\ 
\vdots& \vdots&\vdots&\ddots&\vdots \\
0 & \theta_1^{M-1} & {M-1\choose1}\theta_1^{M-2} & \dots &1  \\ 
\vdots& \vdots&\vdots&\ddots&\vdots   \\
0 & \theta_1^{K-1} & {K-1\choose1}\theta_1^{K-2} & \dots &{K-1\choose  M-1}\theta_1^{K-M} \end{array}
\!\!\right).
\]
In $\widetilde V_n$, i.e.~the block of $M+1$ columns corresponding to $\theta_n$, 
$n = 2, 3,\ldots N$, the first column becomes
$(0\ 1\ \theta_n \ \dots \ \theta_n^{K-1})^T$.  
For the other columns
we perform the additional step, which we used in the initial reduction, i.e.~Gaussian elimination. In this case, we consecutively subtract the previous column and eliminate the common factors. Thus, the element ${k\choose \ell}\theta_n^{k-\ell}$ becomes
\[
\bs
\frac1{\theta_n-\theta_1}&\left({k\choose \ell}\theta_n^{k-\ell} - {k-1\choose \ell}\theta_1\theta_n^{k-1-\ell} - {k-1\choose \ell-1}\theta_n^{k-\ell}\right)\\ & = {k-1\choose \ell}\theta_n^{k-1-\ell}.
\end{split}
\]
Repeating the Vandermonde and Gaussian reductions for each $n=2,3,\ldots, N$ completes the inductive step for $M$.
\epf

\subsection{Sparse Dynamical Sampling}\label{SDS}

Here, we exhibit an example that generalizes on-grid Prony's method and provides a different point of view on some of  the results in \cite{AK16}. For simplicity of exposition, we present it in the case of a finite-dimensional module $\X$.

We consider an abstract initial value problem
\begin{equation}\label{IVP}
	\begin{cases}
	\dot{x}(t)=\A x(t)\\
	x(0)=x_0,
	\end{cases}
	\quad t\in\mathbb R_+,\ x\in\X,
\end{equation} 
where $\X$ is a high-dimensional Euclidean space and $\A \in B(\X)$ is an operator with the 
spectrum $\s(\A)$ and a known basis of generalized eigenvectors $\mathscr B_\A = \{x_\l^m: \l\in\sigma(\A), m = 1, \ldots, m_\l\}$. We consider the standard ordering of the basis so that $\A x_\l^1 = \l x_\l^1$
and $(\A-\l I)x_\l^m = x_\l^{m-1},$ $ m = 2, \ldots m_\l$.

We assume that the initial value $x_0$ is of the form \eqref{signal}, i.e.
\[
x_0 = \sum_{\l\in F}\sum_{m=1}^{m_\l} c_{\l m}x_\l^m,
\]
where $F$ is a relatively small subset of $\s(\A)$. We wish to identify $x_0$ from the dynamical samples
\[
y_\ell(s) = \la x(\beta \ell), e_s\ra, \quad s\in\I, \beta > 0, \ell = 0,1,\ldots L,
\]
where $x(t)$ is the solution of the IVP \eqref{IVP} and $\{e_1, \ldots, e_d\}$ is an orthonormal basis of $\X$. The set $\I\subseteq \{1, \ldots d\}$ is assumed to be such that no eigenvector in $\mathscr B_\A$ is orthogonal to the span of $\{e_s: s\in\I\}$.

\begin{rem}
To recover the classical on-grid Prony problem from this setting, we simply let 
$\A$ be the diagonal operator with the spectrum $\frac{2 \pi i}{d}\cdot\{0,\ldots, d-1\}$, $\mathscr B_\A$ -- the standard  basis, $\{e_1, \ldots, e_d\}$ -- the Fourier  basis, and $\mathcal{I}$ -- a singleton.
\end{rem}

The above problem fits in the framework of dynamical sampling introduced in \cite{ADK13} and further developed in e.g. \cite{ACMT17}. Our goal here is to show that this problem can be solved by Algorithm \ref{alg:rec}. The goal is easily achieved by
putting this problem in the notation of Section \ref{set}. To this end, we let $A: \X\to \CC^{|\I|}$ be defined by $Ax = (\la x, e_s\ra)_{s\in\I}$, and $B\in B(\X)$ -- by
$B = e^{\beta\A}$, so that $x(\beta \ell) = B^\ell x_0$. The function $h: \RR\to \CC$ in \eqref{funh} is then $h(\l) = e^{\beta \l}$ and the representation $\T: \RR\to B(\X)$ is $\T(t) = e^{t\A}$. The non-orthogonality assumption on $\I$ yields that the assumption of Theorem \ref{hardthm} holds and we conclude that Algorithm \ref{alg:rec} can indeed be used to recover all $x_0$ fitting the dynamical samples as long as $\beta$ was chosen in such a way that $h$ is one-to-one on $\s(\A)$. The criteria for uniqueness of such $x_0$ can be derived from \cite{ACMT17}.

\begin{rem}
Strictly speaking we reached the above conclusion prematurely because the group $\T$ will grow exponentially (unless $\s(\A)\subset i\RR$). This violates the non-quasi-analyticity condition \eqref{bdc} and prevents us from using Theorem \ref{spec1}. For this case, however, we can use a much more basic spectral theorem instead, such as \cite[Ch. I, Lemma 3.13]{EN00}. We also note that a similar example can be considered in an infinite dimensional space but it will not provide new insights, beyond those in the previous subsection. Finally, we point out that in this example  Algorithm \ref{alg:rec} does  not need to have access to the spectrum $\s(\A)$ but only to the basis of generalized eigenvectors $\mathscr B_\A$ (provided that $\|\A\|$ is known and $\beta$ was chosen to keep $h$ one-to-one in $\{\l\in\CC: |\l|\le \|\A\|\}$).
\end{rem}

 \section{Identification of time-varying channels}\label{TFS}

 The identification of time-varying channels is a central task in communications engineering.  In mobile communications, for example, the transmission channel is a superposition of various signal paths,  each of which is characterized by a gain factor, a  time-delay and a frequency shift, the latter being a consequence of the Doppler effect. Below, we give a application of the method developed above to identify a transmission channel as described above, that is, the parameters of each signal path.  Clearly, the identification of a communications channel is of utmost importance in order to  transmit information close to channel capacity. For background on the problem see
 \cite{GP14,HeckelBoel13,KP06,PZ14b,PW06,Vtr71,WPK15}.  

In order to apply Prony's method for Banach modules in this setting, we give a formal definition of the channel as follows.   First, we define the time-frequency shift operators $\pi_\l = \pi_{(t,\nu)} = M_{\nu}T_t$   acting on a signal $u:\RR\to\RR$ by 
 $(\pi_\l u)(r) = e^{2\pi i r\nu}u(r+t)$ with $\l =(t,\nu)\in\RR^2$. Observe that $(\pi_\l)^{-1} = e^{2\pi i t\nu}\pi_{-\l}$ \cite{gro01}. 
Second, a generic time-varying communications channel is expressed as the operator
$$ u\mapsto \sum_{\gamma \in F}c_\gamma \pi(\gamma)u\textbf{},$$
where the  gain factor $c_\gamma$ represents the energy transfer on the path (or the collection of paths) with time-frequency shift $\gamma$.
The set $S$ of all potentially appearing time-frequency shifts is a bounded subset of the time-frequency plane, the so called spreading domain. The assumption of boundedness is justified due to physical limitations to the maximal size of time and frequency shifts that a transmission signal may be subjected to. In the current context, the size of the spreading domain is irrelevant and we choose, without loss of generality, the set $S=[0,1)^2$.
Channel identification in this setting, therefore, reduces to the identification of the operator $X = \sum\limits_{\g\in F} c_{\g}\pi_{\g}$, where $F$ is a finite subset in $[0,1)^2$.  A standard method applied to the identification problem is to transmit a well chosen signal $u$ and to recover the unknown parameters of $F$ from the observed channel output $Xu$.

To put this problem into our general framework, we choose $\G = \RR^2$, and, hence,  $\Gg \simeq \RR^2$. We then treat $F\subset [0,1)^2$ as a subset of $\Gg$, and, in order to be compatible with the time-frequency theory, we use the following isomorphism of  $\RR^2$ and $\Gg$: 
$(t_\g, \nu_\g)\mapsto \g \in\Gg$, where
 \begin{equation}\label{griso}
  \g(g)=\g(x,\xi) = e^{2\pi i ( x\nu_\g  - t_\g\xi )},\quad  g = (x,\xi)\in\RR^2.    
 \end{equation}
 In the following, $\HH$ is the Hilbert space $L^2(\RR)$ and $B(\HH)$ is the algebra of all bounded linear operators in $\HH$. 
The space $\X$ is then chosen to be the closure of ${\rm{span}}\{\pi_\l: \l\in\RR^2\}$   in   $B(\HH)$. It is endowed with a Banach module structure by means of the representation $\T: \RR^2\to B(\X)$ given by  $\T(g)X = \pi_g X \pi_{g}^{-1}$, $g\in\RR^2$. 
 
Observe that for $X = \sum\limits_{\g\in F} c_{\g}\pi_{\g}$ with $F\subset [0,1)^2$ finite, we have
 \[
 \T(g)X = \sum_{\g\in F} c_{\g}\pi_g\pi_{\g} (\pi_g)^{-1}= \sum_{\g\in F} c_\g e^{2\pi i (t_{\g}\nu_g - t_g\nu_{\g})}\pi_{\g} = \sum_{\g\in F} c_\g \g(g)\pi_{\g},
 \]
 i.e.,~each $\pi_\g$ is an eigenvector of $\T(g)$ that corresponds to the eigenvalue $\g(g)$. This justifies our choice of isomorphism in \eqref{griso}. It also follows that $\L(\pi_\g,\T)=\{\g\}$. Moreover, the spectral submodules $\X_\g$ are one-dimensional. Thus, we showed that time-varying communication channels are a special case of signals of the form \eqref{signal}. 

To define the operator $A: \X\to \CC^S$, we fix the function $u\in \HH$ to be the Gaussian
\[
u(t) = e^{-t^2}, \ t\in \RR,
\]
and pick a finite subset $K\subset \RR^2$ of cardinality $S$. We let
\[
(AX)(s) = \la X\pi_s u, \pi_s u\ra,\ s\in K.
\]

We note that similar types of measurements were used in \cite{GP14}, but the recovery algorithm developed there is different.

For $B$ of the form \eqref{opB}, we have $h(\g) = \sum\limits_{n=1}^Nb_n\g(g_n)$, as before.
We can still choose $B = \T(g_1)$ with $g_1 = (1,1)$. In this case, however, the function $h(\g) = \g(g)$ is not one-to-one on $\Omega = [0,1)^2$. Nevertheless, the measurements
\[
y_\ell (s) = (AB^\ell X)(s)  = \sum_{n=1}^N\sum_{\g\in F} b_nc_\g \g^\ell(g_n) \la \pi_\g\pi_s u, \pi_s u\ra
\]
still allow one to find the polynomial $p_{\min}$ with the set $R_{\min}$ of non-zero roots guaranteed to be a subset of $h(F)$. Observe also that with this choice of $A$ and $B$ all measurements $y_\ell(s)$ are linear combinations of a finite collection of numbers of the form $\la X\pi_r u, \pi_q u\ra$, for some $r,q\in\RR^2$.

Restricting the set $\Omega$ in such a way that $h$ satisfies the assumptions of Theorem \ref{hardthm} yields $R_{\min} = h(F)$. We also remark that 
choosing several different $B$'s one can get $F$ after  recovering the corresponding subsets of $h(F)$.

\begin{exmp}
Choosing $N=2$, $b_1=1$, $b_2=i$, $g_1=(0,1/12)$, and
$g_2=(1/12,0)$, we get the function
\begin{equation}\label{goodh}
    h(t,\nu)=e^{\frac{2\pi i t}{12}}+i e^{-\frac{2\pi i \nu}{12}}= \cos\frac{\pi t}6+i\sin\frac{\pi t}6 +i \cos\frac{\pi \nu}6  +  \sin\frac{\pi \nu}6
\end{equation}
which is injective on $[0,1)^2$.  To see that, assume that $t, t', \nu, \nu'\in [0,1)$ are such that $h(t,\nu) = h(t',\nu')$. Then we must have
\[
\begin{cases}
 \cos(\pi t/6)+  \sin(\pi \nu/6) = \cos(\pi t'/6)+  \sin(\pi \nu'/6) \\
\sin(\pi t/6)+  \cos(\pi \nu/6) = \sin(\pi t'/6)+  \cos(\pi \nu'/6)
\end{cases}
\]
or, equivalently,
\[
\begin{cases}
 \cos(\pi t/6)- \cos(\pi t'/6)=  \sin(\pi \nu'/6) - \sin(\pi \nu/6) \\
\sin(\pi t/6)-\sin(\pi t'/6)   =   \cos(\pi \nu'/6)-\cos(\pi \nu/6).
\end{cases}
\]
Using the standard trig identities the above is equivalent to
\[
\begin{cases}
 \sin\frac{\pi(t+t')}{12}\sin\frac{\pi(t-t')}{12}=  \cos\frac{\pi(\nu+\nu')}{12}\sin\frac{\pi(\nu-\nu')}{12} \\
 \cos\frac{\pi(t+t')}{12}\sin\frac{\pi(t-t')}{12}=  \sin\frac{\pi(\nu+\nu')}{12}\sin\frac{\pi(\nu-\nu')}{12}.
\end{cases}
\]
Observe that for $t, t', \nu, \nu'\in [0,1)$ both sides of the last equation can only equal $0$ if $t=t'$ and $\nu=\nu'$, which would not contradict the injectivity of $h$. Otherwise, we can divide the first equation by the second one and obtain
\[
\tan\frac{\pi(t+t')}{12} = \cot \frac{\pi(\nu+\nu')}{12},
\]
which would only hold if
\[
\frac{\pi(t+t')}{12} = \frac{\pi}{2} - \frac{\pi(\nu+\nu')}{12} +\pi k
\]
for some $k\in\ZZ$. But that would require $t+t'+\nu+\nu' = 6+12k$ which is impossible for $t, t', \nu, \nu'\in [0,1)$. We conclude that the function $h$ in \eqref{goodh} is, indeed, injective on $[0,1)^2$. In fact, an elementary computation yields that the inverse function $h^{-1}: h([0,1)^2)\subset\CC\to [0,1)^2$ is given by
\begeq
h^{-1}(x+iy) = \frac{3}{\pi}
\begin{pmatrix}
\arcsin{\frac{x^2+y^2-2}{2}}-\arcsin{\frac{x^2-y^2}{x^2+y^2}}\\
\arcsin{\frac{x^2+y^2-2}{2}}+\arcsin{\frac{x^2-y^2}{x^2+y^2}}
\end{pmatrix}.
\eq

\end{exmp}

To describe a more general model, we let $\HH_n$ be the $n$-th order Sobolev space and define the following two operators. We let $D:  \HH_1\to \HH$ be the derivative operator: $(Du)(t) = \dot u(t)$ and $\Delta: \F\HH_1\to \HH$ be given by $\Delta = \F^{-1}D\F$, where $\F\in B(\HH)$ is the classical unitary Fourier transform. The space $\X$ will then be a subspace of $L(\HH_M\cap\F\HH_M, \HH)$ spanned by $\Delta^j D^k\pi_\g$, $\g\in\RR^2$, $j, k \in\{0,1,\ldots, M\}$. The representation $\T$ remains to be defined in the same way, albeit with a different range space. Since
\[
\T(g)\Delta = \pi_g \F^{-1}D\F(\pi_g)^{-1}= 2\pi i t I +\Delta\ \mbox{ and}
\]
\[
\T(g)D = \pi_g D(\pi_g)^{-1}= 2\pi i \nu I +D,\ g = (t,\nu),
\]
we have $\L(D,\T) = \L(\Delta,\T) = \{0\}$, and it follows \cite[Corollary 7.8]{BK05} that $\L(\Delta^jD^k\pi_\g) = \{\g\}$, $\g\in\RR^2$, $j, k \in\{0,1,\ldots, M\}$.

Thus, in this setting, one seeks to recover an operator of the form
\[
X = \sum_{\g\in F}q^{\Delta}_{\g}(\Delta)q^D_{\g}(D)\pi_\g, 
\]
where $q^{\Delta}_\g$ and $q^D_\g$ are non-zero polynomials of degree at most $M$.
The same choice of $A$ and $B$ as earlier in this section, once again, allows one to find a polynomial $p_{\min}$ whose non-zero roots are guaranteed to be in $h(F)$, where $h$ is still given by \eqref{funh}. As before, this is a direct consequence of Proposition \ref{mainprop} and Theorem \ref{spec1}.

\section{Concluding remarks}

In this short note, we presented only a few instances of how our point of view on Prony's method can be applied.
Many more remain to be developed.  
One problem that can be investigated is identifying rational functions using the fact that they are Laplace transforms of almost periodic functions with a finite spectrum. Another problem is identifying polynomial splines using the  observation that their high order derivatives are linear combinations of Dirac masses. We invite the readers to pursue these and other related problems.

\newpage
\centerline{\bf Acknowledgements}
\medskip
This paper is dedicated to Karlheinz Gr\"ochenig. Charly has not only been the leader of the vast area of time-frequency analysis since the 1990s, but also significantly influenced the work of the authors of this paper. We have  greatly appreciated his direct, honest, and unapologetic criticism
as well as the substantial (albeit rare) praise 
for the handful of our papers that actually deserved to be published (in Charly's eyes). Above all, we value our ongoing friendship, with hopes that it will endure for many more years to come.

The first author acknowledges the support by the NSF grant DMS-2208031 and thanks the second author for his hospitality during the visits to KU Eichst\"att-Ingolstadt. The second author acknowledges support by the DFG grants PF 450/11-1 and PF 450/9-1.

\bibliographystyle{siam}
\bibliography{refs}

\end {document}